\documentclass[11pt,reqno,oneside]{amsart}
 \usepackage{verbatim,amsmath,amssymb,cite,xspace}
\usepackage{color,cite,graphicx}

\numberwithin{equation}{section}
\DeclareMathOperator*{\esssup}{ess\,sup}
\theoremstyle{plain}
\newtheorem{theorem}{Theorem}[section]
\newtheorem{lemma}{Lemma}[section]

\theoremstyle{definition}
\newtheorem{definition}{Definition}[section]

\theoremstyle{remark}

\begin{document}
\title[]{Are the incompressible 3d Navier-Stokes equations  locally ill-posed in the natural energy space?}
\author{Hao Jia\,,\, Vladimir Sverak \\ \\
\smaller University of Minnesota}
\maketitle
 \begin{abstract}{An important open problem in the theory of the Navier-Stokes equations is the uniqueness of the Leray-Hopf weak solutions with $L^2$ initial data.
 In this paper we give sufficient conditions for non-uniqueness
 in terms of spectral properties of a natural linear operator associated to scale-invariant solutions recently constructed in \cite{JiaSverak}. If the spectral conditions are satisfied,  non-uniqueness and ill-posedness can appear for quite benign compactly supported data, just at the borderline of applicability of the classical perturbation theory.  The verification of the spectral conditions seems to be approachable by relatively straightforward numerical simulations which involve only smooth functions.}
  \end{abstract}
\begin{section}{Introduction}
We consider the classical Cauchy problem for the 3d Navier-Stokes equations (NSE) in $R^3\times (0,\infty)$,
\begin{eqnarray}
\begin{array}{rl}
\partial_tu-\Delta u+u\cdot \nabla u+\nabla p&=0\,,\\
{\rm div~}u&=0\,,\\
u(\cdot,0)&=u_0\,.
\end{array}
\end{eqnarray}
Let us recall the scaling symmetry
\begin{eqnarray*}
u(x,t)&\rightarrow&u_{\lambda}(x,t):=\lambda u(\lambda x,\lambda t),\\
u_0(x)&\rightarrow& u_{0\lambda}(x):=\lambda u_0(\lambda x),\\
p(x,t)&\rightarrow&p_{\lambda}(x,t):=\lambda^2p(\lambda x,\lambda^2 t).
\end{eqnarray*}
defined for $\lambda>0$. If $u,p,u_0$ satisfy the equations,  so do $u_\lambda, p_\lambda, u_{0\lambda}$.
 A solution $u$ is {\it scale invariant} if $u\equiv u_{\lambda}$ for all $\lambda>0$. The initial condition $u_0$ is scale invariant if $u_0\equiv u_{0\lambda}$ for all $\lambda>0$. In \cite{JiaSverak}, it was proved that for each scale invariant initial data $u_0\in C^{\alpha}(R^3\backslash \{0\})$ there exists at least one scale invariant solution $u\in C^{\infty}(R^3\times (0,\infty))\cap C^{\alpha}(R^3\times [0,\infty)\backslash\{(0,0)\})$. See also \cite{Tsai} for a generalization to discretely self-similar solutions. 
 
 It has been conjectured in~\cite{JiaSverak, HAOSVERAK} that for many large scale-invariant initial data the scale invariant solutions are not unique, and possible implications for the non-uniqueness of the Leray-Hopf weak solutions were suggested.
  In this paper we investigate these topics further.
  
 Let us consider a scale-invariant initial condition $u_0$ which is smooth away from the origin. For $\sigma\ge 0$, at first taken sufficiently small so that we have uniqueness, 
   let $u_{\sigma}(x,t)=\frac{1}{\sqrt{t}}U_{\sigma}(\frac{x}{\sqrt{t}})$ be the scale invariant solution to NSE with the initial data $\sigma u_0$. The field $U_{\sigma}$ satisfies
\begin{equation}\label{eq:mainselfsimilarequation}
\Delta U_{\sigma}+\frac{x}{2}\cdot \nabla U_{\sigma}+\frac{1}{2}U_{\sigma}-U_{\sigma}\cdot\nabla U_{\sigma}+\nabla P=0
\end{equation}
in $R^3$, with $|U_{\sigma}(x)-\sigma u_0(x)|=o(\frac{1}{|x|})$ as $x\to \infty$. By the results in \cite{JiaSverak} we know $|U_{\sigma}(x)-u_0(x)|=O(\frac{1}{|x|^3})$ as $|x|\to \infty$.
We consider Navier-Stokes solutions $u(x,t)$ which are close to $u_\sigma(x,t)$ and  write them as
\begin{equation}\label{vs1}
u(x,t)={1\over \sqrt t}\, U_\sigma\left({x\over\sqrt t}\right)+
{1\over \sqrt{t}}\,\phi\left(\,{x\over\sqrt t}\,,\,t\,\right)\,.
\end{equation}
The equation for $\phi$ is
\begin{equation}\label{vs2}
t\phi_t=\Delta\phi+{x\over 2}\cdot\nabla\phi+{1\over 2}\phi-U_\sigma\cdot\nabla\phi-\phi\cdot\nabla U_\sigma-\phi\cdot\nabla\phi+\nabla \pi\,,
\end{equation}
for a suitable  function $\pi$ (related to the pressure).

The linearization of this equation will be written as
\begin{equation}\label{vs3}
t\phi_t=\mathcal{L}_\sigma\phi\,,
\end{equation}
where
\begin{equation}\label{vs20}
\mathcal{L}_{\sigma}\phi:=\Delta \phi+\frac{x}{2}\cdot\nabla\phi+\frac{1}{2}\phi-U_{\sigma}\cdot\nabla \phi-\phi\cdot\nabla U_{\sigma}+\nabla P\,,
\end{equation}
 and the function $P$ is chosen so that $\mathcal{L}_\sigma\phi$ is divergence-free.  Also, $P$ is assumed to have a suitable decay at $\infty$ (so that it is uniquely determined, perhaps up to a constant).
 We expect that for $\phi$ which is small in suitable Banach spaces  the behavior of the solutions~(\ref{vs1}) is determined by the linearized equation~(\ref{vs3}), which in turn should be governed by the spectral properties of $\mathcal{L}_\sigma.$
 \newcommand{\lsig}{\mathcal{L}_\sigma}
  At a heuristic level, if $\Phi$ is an eigenvector of $\lsig$ with an eigenvalue $\lambda$ which has a strictly positive real part and $\Phi(x)=o(\,|x|^{-1})$ for $|x|\to\infty$, then~(\ref{vs3}) $t^\lambda\,\Phi$ is a solution of~(\ref{vs3}) with initial value $0$ indicating non-uniqueness and ill-posedness.
  
 In what follows we will often replace the time $t$ in~(\ref{vs2}) with the non-dimensionalised time $\log {t\over t_0}$ for some reference time $t_0$, or simply $\log t$, with slight abuse of the usual rules for handling quantities which are not dimensionless. With the new time $t$ equation~(\ref{vs3}) becomes
 \begin{equation}\label{vs3b}
 \phi_t=\lsig\phi, \qquad t\in (-\infty, 0)\,.
 \end{equation}

As we work in the whole space, we can expect continuous spectra and the dependence of spectral properties of $\lsig$ on the  Banach spaces in which the problem is considered.

The choice of the spaces is dictated by the ``boundary conditions" of $\phi$ near $x\sim\infty$. The conditions for $U_\sigma$ are (see~\cite{JiaSverak})
\newcommand{\be}{\begin{equation}}
\newcommand{\ee}{\end{equation}}
\newcommand{\bb}{\bigskip}
\newcommand{\la}{\label}
\be\la{vs5}
U_\sigma(x)=\sigma u_0(x)+O\left({1\over\,\,\,|x|^{3}}\right)\,, \qquad |x|\to\infty.
\ee 
When studying possible non-uniqueness of the solution of the steady equation~(\ref{eq:mainselfsimilarequation})  we should impose on the perturbations a condition consistent with
\be\la{vs6}
\phi(x)=O\left({1\over\,\,\,|x|^{3}}\right)\,,\qquad |x|\to\infty\,.
\ee
 For the more general time-dependent perturbations $\phi(x,t)$ the most general condition which still gives non-uniqueness of the initial-value problem is
\be\la{vs7}
\phi(x,t)=o\left({1\over |x|}\right)\,,\qquad |x|\to\infty\,.
\ee

In this paper we will work with the space
\begin{equation}
X:=\{\phi\in L^2\cap L^4(R^3):~{\rm div~}\phi=0\},
\end{equation}
with the natural norm
\begin{equation}
\|\phi\|_X:=\|\phi\|_{L^2(R^3)}+\|\phi\|_{L^4(R^3)}.
\end{equation}
\newcommand{\rf}[1]{(\ref{#1})}
This accommodates condition~\rf{vs6} although not quite the condition~\rf{vs7}. However, it is not hard to see that our arguments work, {\it mutatis mutandis}, if we replace $X$ by $L^{3-\delta}\cap L^{3+\delta}$ (together with the divergence free condition) for $\delta>0$. The case $\delta=0$ might require additional effort.

Recalling the definition of $\lsig$ in~\rf{vs20}, we let ${\mathcal L}={\mathcal L}_0$.
We define the domain of $\mathcal{L_{\sigma}}$ as
\begin{equation}
\mathcal{D}:=\{\phi\in X: \partial_j\phi, \partial_{ij}\phi\in X, {\rm ~and~} x\cdot\nabla\phi\in X\}.
\end{equation}
 We will show below that the spectrum of $\mathcal{L}$ is contained in the set
\begin{equation}
\Sigma:=\{\lambda\in\mathbb{C}:{\rm ~Re\,}\lambda\leq -\frac{1}{4}\}
\end{equation}
and $K_{\sigma}\phi:=U_{\sigma}\cdot\nabla \phi+\phi\cdot\nabla U_{\sigma}+\nabla P$ is relatively compact with respect to $\mathcal{L}$. Thus the spectrum of $\mathcal{L}_{\sigma}=\mathcal{L}-\mathcal{K}(\sigma)$ is the union of
one part contained in $\Sigma$ and some isolated eigenvalues in the set
\begin{equation*}
\{\lambda\in\mathbb{C}:-\frac{1}{4}<{\rm Re\,}\lambda\}.
\end{equation*}
When $\sigma$ is small, one can solve equation (\ref{eq:mainselfsimilarequation}) by perturbation argument (and the solution is unique). Moreover, still for small $\sigma$, the set of eigenvalues of $\mathcal{L}_{\sigma}$ can be shown to remain away from the imaginary axis, with ${\rm Re}<0$. It is not hard to show that as long as $0$ is not an eigenvalue of $\mathcal{L}_{\sigma}$, we can use perturbation arguments to continue the solution curve $U_{\sigma}$ as a regular function of  $\sigma$. If we increase $\sigma$ further, some eigenvalues of $\mathcal{L}_{\sigma}$ might cross the imaginary axis. We consider the following scenarios for the crossing (which cover the ``generic case"):

\medskip
({\bf A}) For $\sigma\leq\sigma_0$, the eigenvalues of $\mathcal{L}_{\sigma}\subseteq\{\lambda\in\mathbb{C}:{\rm ~Re}\,\lambda<-\delta<0\}\cup\{\lambda_1(\sigma),~\overline{\lambda_1(\sigma)}\}$ for some $\delta>0$ and ${\rm Re}\,\lambda_1(\sigma)< 0$ for $\sigma<\sigma_0$. $\lambda_1(\sigma),\,\overline{\lambda_1(\sigma)}$ are simple eigenvalues. Moreover, ${\rm Re}\,\lambda_1(\sigma_0)=0$, ${\rm Im}\,\lambda_1(\sigma_0)>0$ and $\frac{d}{d\sigma}|_{\sigma_0}{\rm Re}\,\lambda_1>0$.

\medskip

or

\medskip
({\bf B}) For $\sigma\leq\sigma_0$, the eigenvalues of $\mathcal{L}_{\sigma}\subseteq\{\lambda\in\mathbb{C}:{\rm ~Re}\,\lambda<-\delta<0\}\cup\{\lambda_1(\sigma)\}$ for some $\delta>0$ and $\lambda_1(\sigma)< 0$ for $\sigma<\sigma_0$. $\lambda_1(\sigma)$ is a simple eigenvalue. Moreover, $\lambda_1(\sigma_0)=0$.

\medskip

Under assumption ({\bf A}) we show that we can continue the solution $U_{\sigma}$ as a regular curve to some $\sigma_1>\sigma$. In addition, we have the following theorem.
\begin{theorem}\label{th1}
Suppose ({\bf A}) holds.  Let $\sigma_1>\sigma_0$ be sufficiently close to $\sigma_0$. Then for $\sigma_1>\sigma>\sigma_0$ and any small $\epsilon>0$, there exists an ancient solution $\phi=\phi(x,t,\sigma)$ to
\begin{equation}
\left\{\begin{array}{rl}
\partial_t\phi&=\Delta\phi+\frac{x}{2}\cdot\nabla\phi+\frac{\phi}{2}-\phi\cdot\nabla U_{\sigma}-U_{\sigma}\cdot\nabla\phi-\phi\cdot\nabla\phi+\nabla P\\
{\rm div~}\phi&=0\\
         \end{array}\right.~{\rm in~} R^3\times(-\infty,0).
\end{equation}
with
\begin{equation}
\|\phi(\cdot,t,\sigma)\|_X+ \|\nabla\phi(\cdot,t,\sigma)\|_X<\epsilon, {\rm ~for~}t\in(-\infty,0).
\end{equation}
Moreover,
\begin{equation}
\inf_{t\in (-\infty,0)}e^{-\frac{1}{16}t}\|\phi(\cdot,t,\sigma)\|_{L^4(R^3)}>0.
\end{equation}
\end{theorem}

\smallskip
\noindent
Once we have such solutions, we define
\begin{equation}
u(x,t):=\frac{1}{\sqrt{t}}\phi(\frac{x}{\sqrt{t}},\log{t},\sigma)+\frac{1}{\sqrt{t}}U_{\sigma}(\frac{x}{\sqrt{t}}), {\rm ~for~}x\in R^3,~t>0.
\end{equation}
The field $u$ satisfies NSE with initial data $\sigma u_0$ and
\begin{equation}
\liminf_{t\to 0+}\|u(\cdot,t)-\frac{1}{\sqrt{t}}U_{\sigma}(\frac{x}{\sqrt{t}})\|_{L^4(R^3)}\ge \liminf_{t\to 0+}t^{-\frac{1}{8}}\|\phi(\cdot,\log{t})\|_{L^4(R^3)}=+\infty.
\end{equation}
\\
Under assumption ({\bf B}) and some additional generic non-degeneracy conditions (see Section 4 for details), we show there is a ``saddle-node" bifurcation at $\sigma=\sigma_0$, which gives another self similar solution $\frac{1}{\sqrt{t}}\tilde{U}_{\sigma}(\frac{x}{\sqrt{t}})$, for $\sigma<\sigma_0$ sufficiently close to $\sigma_0$.\\

In summary, under the spectral assumption ({\bf A}) or ({\bf B}) we have two solutions to Navier Stokes equation with initial data $\sigma u_0$ when $\sigma$ is close to $\sigma_0$ and $\sigma>\sigma_0$ (in the case of ({\bf A})), or $\sigma<\sigma_0$ (in the case of ({\bf B})). 

\bb
The solutions obtained above have infinite energy, due to the slow decay at $\infty$. Our next goal is to localize them to obtain non-uniqueness for Leray-Hopf  solutions. The reason why this is not straightforward is  that the equations for the perturbations contain critically singular lower order terms which cannot be assumed to be small.  The following model problem illustrates the difficulty which has to be overcome. We consider the linear heat equation with a potential
\newcommand{\dd}{\Delta}
\be\la{vs8}
u_t-\dd u - V(x,t)\,u=0\,. 
\ee

The classical parabolic techniques show that the equation can be considered as a perturbation of the case $V=0$ under the assumption 
\be\la{vs9}
V\in L_{t,x}^{n+2\over 2}\,.
\ee
In this case equation~\rf{vs9} is well-posed for initial data in $L^p$ for any $p\in (1,\infty)$, by usual perturbation arguments.
We note that the scaling
\be\la{vs9}
V(x,t)\to\lambda^2 V(\lambda x,\lambda^2 t)\,
\ee
preserves the norm in $L_{t,x}^{n+2\over 2}$. Let us now assume that $V$ is of the form
\be\la{vs11}
V(x,t)={1\over t}A\left(x\over \sqrt t\right)
\ee
for some smooth function $A$ on $R^n$ with suitable decay at $\infty$. Then $V$ is invariant under the scaling~\rf{vs11} and ``just misses" the critical space
$L_{t,x}^{n+2\over 2}$. The well-posedness (in $L^p$, say) is then not clear from the general parabolic theory. One can hope to replace $L_{t,x}^{n+2\over 2}$ with some slightly weaker critical space $Z$ which contains~\rf{vs11}, and this may indeed work  if $||V||_Z$ is sufficiently small. Such results are in fact essentially optimal, as~\rf{vs8} with $V$ given by~\rf{vs11} may in general be ill-posed in $L^p$. We can seek the solutions of~\rf{vs8} in the form
\be\la{vs12}
u(x,t)={1\over t^{\alpha\over 2}}\,\phi\left({x\over \sqrt t},t\right)\,,
\ee
which gives
\be\la{vs13}
t\phi_t-{1\over 2}x\cdot\nabla \phi-{\alpha\over 2}\phi-\dd\phi-A\phi=0\,.
\ee
Heuristically we expect that a necessary condition for well-posedness of this equation is that the the spectrum of the operator
\be\la{vs14}
\phi\to\dd\phi+{1\over 2}x\cdot\nabla\phi+A\phi+{\alpha\over 2}\phi
\ee
be strictly contained in the left half-plane 
\be\la{vs15}
\{z\,\,,\,\,\Re z<0\}\,.
\ee
This is the same as the spectrum of 
\be\la{vs16}
\phi\to\dd\phi+{1\over 2}x\nabla \phi+A\phi
\ee
being contained in
\be\la{vs16b}
\{z,\,\,\Re z<-{\alpha\over 2}\}\,.
\ee
We note that 
\be\la{vs17}
||{1\over t^{\alpha\over 2}}
\phi\left({\,\cdot\,\over \sqrt t}\right)||_
{L^p(R^n)}
=t^{{n\over 2p}-{\alpha\over 2}}||\phi||_{L^p(R^n)}\,,
\ee
and therefore an eigenvalue $-{\alpha\over 2}$ (with suitable decay of the eigenfunction at $x\sim\infty$) suggests the local ill-posedness of the Cauchy problem in $L^p(R^n)$ with $p<{n\over\alpha}\,$. On the other hand if the spectrum is contained in~\rf{vs16b}, we can expect the local well-posedness in $L^p(R^n)$ for $p>{n\over\alpha}$. We see that there is an interplay between the spectrum of~\rf{vs16} and the spaces in which equation~\rf{vs8} is well-posed.
 
We will apply these ideas in the context of equation~\rf{vs3}.
Under the spectral assumption ({\bf A}) or ({\bf B}) this program can be carried out.

We introduce the space
\begin{equation}
Y:=\{a\in L^{\infty}(R^3)|{\rm ~div~}a=0, {\rm ~and~}\esssup_x\left((1+|x|)^{1+|\alpha|}|\nabla^{\alpha}a(x)|\right)<\infty, {\rm ~for~}|\alpha|=0,1,2.\}.
\end{equation}
equipped with the natural norm
\begin{equation}
\|a\|_Y:=\esssup_{x\in R^3,|\alpha|=0,1,2}\left((1+|x|)^{1+|\alpha|}|\nabla^{\alpha}a(x)|\right).
\end{equation}
For $a\in Y$, define operator
\begin{equation}
\mathcal{K}(a)\phi:=a\cdot\nabla \phi+\phi\cdot\nabla a+\nabla p, {\rm ~for~any~}\phi\in \mathcal{D}.
\end{equation}
\begin{theorem}
Let $a\in Y$ be such that $\frac{1}{8}>\beta:=s(\mathcal{L}-\mathcal{K}(a))$ which is the maximum of real part of any eigenvalue of $\mathcal{L}-\mathcal{K}(a)$. Denote $\tilde{a}(x,t)=\frac{1}{\sqrt{t}}a(\frac{x}{\sqrt{t}})$. Let
\begin{equation}
\tilde{b}(x,t):=\frac{1}{\sqrt{t}}b(\frac{x}{\sqrt{t}},\log{t}), {\rm ~for~}x\in R^3,~t\in (0,1).
\end{equation}
Suppose
\begin{equation}
\epsilon=\esssup_{s\in(-\infty,0)}(\|b(\cdot,s)\|_X+\|\nabla b(\cdot,s)\|_X)
\end{equation}
is sufficiently small depending on $\|a\|_Y,\,\beta$. Let $T=T(\|a\|_Y,\beta,\|u_0\|_{L^4(R^3)})>0$ be sufficiently small.
Then there exists a unique solution $u\in Z_T$ to the generalized Navier Stokes system with singular lower order terms
\begin{eqnarray}\label{eq:singularstokes}
\left.\begin{array}{rl}
    \partial_tu-\Delta u+\tilde{a}\cdot\nabla u+u\cdot\nabla\tilde{a}+\tilde{b}\cdot\nabla u+u\cdot\nabla \tilde{b}+u\cdot\nabla u+\nabla p&=0\\
                                                                  {\rm div~}u&=0\\
         \end{array}\right\} {\rm ~in~}R^3\times(0,T),
\end{eqnarray}
with initial data $u(\cdot,0)=u_0\in L^4(R^3)$ in the sense that
\begin{equation}
\lim_{t\to0+}\|u(\cdot,t)-u_0\|_{L^4(R^3)}=0.
\end{equation}
In the above,
\begin{equation}
Z_T:=\{u\in L_t^{\infty}L_x^4(R^3\times(0,T)): ~\esssup_{t\in (0,T)}t^{\frac{1}{2}}\|\nabla u(\cdot,t)\|_{L^4(R^3)}<\infty, ~\lim_{t\to0+}t^{\frac{1}{2}}\|\nabla u(\cdot,t)\|_{L^4(R^3)}=0\},
\end{equation}
equipped with the natural norm
\begin{equation}
\|u\|_{Z_T}:=\esssup_{t\in (0,T)}\left(\|u(\cdot,t)\|_{L^4(R^3)}+t^{\frac{1}{2}}\|\nabla u(\cdot,t)\|_{L^4(R^3)}\right).
\end{equation}
Moreover, $u$ satisfies
\begin{equation}
\|u\|_{L_t^{\infty}L_x^4(R^3\times(0,T))}+\sup_{t\in(0,T)}t^{\frac{1}{2}}\|\nabla u(\cdot,t)\|_{L_x^4(R^3)}\leq C(\beta,\|a\|_Y)\|u_0\|_{L_x^4(R^3)}.
\end{equation}
\end{theorem}
Using this theorem we can truncate the scale-invariant solutions with initial data $\sigma u_0$, and obtain the following result.
\begin{theorem}
Suppose ({\bf A}), or ({\bf B}) together with non-degeneracy conditions in Theorem \ref{th:bifurcationB} holds. Then there exist two different Leray-Hopf weak solutions $u_1$ and $u_2$ to NSE with compactly supported initial data $v_0$ which is smooth in $R^3\backslash\{0\}$ and near origin
$|v_0(x)|=O(\frac{1}{|x|})$.
\end{theorem}

The question of uniqueness of the Leray-Hopf solutions has a long history.
E.~Hopf's comment on this issue in his 1951 paper~\cite{EH}, p. 217, is: ``It is hard to believe that the initial value problem for the viscous fluid in dimension $n=3$ could have more than one solution, and more work should be devoted to the study of the uniqueness question."\footnote{``Es ist
schwer zu glauben, dass die Anfangs\-wert\-aufgabe z\"aher Flussig\-keiten
f\"ur $n = 3$ mehr als eine L\"osung haben k\"onnte, und der Erledig\-ung
der Eindeutig\-keits\-frage sollte mehr Auf\-merksam\-keit geschenkt
werden."}

O.\ A.\ Ladyzhenskaya wrote in 1969 on the issue of uniqueness, see~\cite{OAL}, p. 229: ``As regards the class of weak Hopf solutions for the general three-dimensional case, it has always seemed to me that it is too broad, i.\ e.\ , that there is missing in it a basic property of of the initial-value problem, viz.\ its determinacy (a uniqueness theorem)... But I had available only indirect reasons in support of this assertion... which had no formal demonstrative power. At this time I am able to rigorously prove the validity of my opinions." She then proceeds to construct an example of non-uniqueness in a non-standard time-dependent domain which degenerates to one point at time $t=0$. The boundary conditions are also non-standard. As our constructions here, Ladyzhenskaya's construction is  based on the study of certain
scale-invariant solutions. The setup is different from ours and produces two distinct solutions with the same boundary conditions and the same non-zero right-hand-side. The solutions are considered in the class of axi-symmetric velocity fields with no swirl.
The time-dependent domain corresponds to the image of a fixed space-time domain in the self-similarity variables. A comment on p. 233 reads: ``We note that a certain `exoticness' of the domain $Q_T$ in which our example has been constructed does not imply the loss of the uniqueness theorems just mentioned (they are usually proved for the case of a domain $\Omega_1$ which does not alter with time\,[\footnote{The uniqueness theorems mentioned here refer to the Ladyzhenskaya-Prodi-Serrin type uniqueness results. }]...The example described here can provoke `displeasure' for only one reason. It has been constructed for boundary conditions (18), but not for adhesion conditions..."
If the spectral conditions we assume in the present paper are satisfied, then all the undesirable features of Ladyzhenskaya's example can be removed and one will have a more or less optimal example of the local-in-time ill-posedness in the energy space. The example will be at the borderline of the class in which uniqueness can be proved via the Ladyzhenskaya-Prodi-Serrin criteria.

\end{section}

\begin{section}{Properties of $\mathcal{L}$}
In this section we study the spectral properties of operator $\mathcal{L}$ and the semigroup it generates. Firstly we prove the following lemma.
\begin{lemma}\label{lm:spectrumofmainoperator}
Let $X, \mathcal{L}$ and $\mathcal{D}\subseteq X$ be defined as above. Then $\mathcal{L}$ is a densely defined, closed operator. Moreover the resolvent set
\begin{equation}
\rho(\mathcal{L})\supseteq \{\lambda\in\mathbb{C}: {\rm Re\,}\lambda >-\frac{1}{4} \}.
\end{equation}
\end{lemma}

\smallskip
\noindent
{\bf Proof.} Let $\phi\in \mathcal{D}$ with $\mathcal{L}\phi-\lambda \phi=f$, that is
\begin{equation*}
\Delta\phi+\frac{x}{2}\cdot\nabla\phi+\frac{1}{2}\phi-\lambda \phi=f,
\end{equation*}
where $\lambda=\beta+i\gamma$ and $\beta>-\frac{1}{4}$.
Introduce
\begin{equation*}
h(x,t):=t^{-\frac{1}{2}+\beta+i\gamma}\phi(\frac{x}{\sqrt{t}}).
\end{equation*}
Since
\begin{equation*}
\|h(\cdot,t)\|_{L^2(R^3)}=t^{\frac{1}{4}+\beta}\|\phi\|_{L^2(R^3)}\to 0, {\rm ~as~}t\to 0,
\end{equation*}
it is easy to check $h$ satisfies
\begin{eqnarray*}
\left.\begin{array}{rl}
\partial_th-\Delta h&=-t^{-\frac{3}{2}+\beta+i\gamma}f(\frac{x}{\sqrt{t}})\\
 u(\cdot,0)&=0
\end{array}\right\} {\rm ~in~}R^3\times(0,\infty).
\end{eqnarray*}
Note that
\begin{equation*}
t^{-\frac{3}{2}+\beta}\|f(\frac{\cdot}{\sqrt{t}})\|_{L^2_x(R^3)}\leq t^{-\frac{3}{4}+\beta}\|f\|_X
\end{equation*}
is integrable in $t$ locally. Thus we can write
\begin{equation*}
h(\cdot,t)=-\int_0^te^{\Delta(t-s)}s^{-\frac{3}{2}+\beta+i\gamma}f(\frac{\cdot}{\sqrt{s}})ds,
\end{equation*}
and consequently
\begin{equation}
\|u(\cdot,\frac{1}{2})\|_{L^2(R^3)}\leq C\int_0^{\frac{1}{2}}s^{-\frac{3}{4}+\beta}\|f\|_{L^2(R^3)}ds\leq C(\beta)\|f\|_{L^2(R^3)}.
\end{equation}
Now for $t\ge \frac{1}{2}$ write $u$ as
\begin{equation}
u(\cdot,t)=e^{\Delta (t-\frac{1}{2})}u(\cdot,\frac{1}{2})-\int_{\frac{1}{2}}^te^{\Delta (t-s)}s^{-\frac{1}{2}+\beta+i\gamma}f(\frac{\cdot}{\sqrt{s}})ds.
\end{equation}
Note that in the above formula the nonhomogeneous term is no longer singular in $s$. We obtain
\begin{eqnarray*}
\|u(\cdot,\frac{3}{4})\|_X&\leq&C\|u(\cdot,\frac{1}{2})\|_{L^2(R^3)}+C(\beta)\|f\|_X\\
&\leq&C(\beta)\|f\|_X.
\end{eqnarray*}
Then it is not difficult to deduce from the general theory of heat equation that
\begin{equation}
\||\partial_tu|+|\nabla^2 u|\|_{L^2_tX_x(R^3\times[\frac{7}{8},1])}\leq C(\beta)\|f\|_X.
\end{equation}
Now note the relation between $u$ and $\phi$, we see
\begin{equation}\label{eq:estimateforlinearoperator}
\|\phi\|_X+\|\nabla\phi\|_X+\|\nabla^2\phi\|_X+\|\frac{x}{2}\cdot\nabla\phi-i\gamma\phi\|_X\leq C(\beta)\|f\|_X~{\rm for~}\beta>-\frac{1}{4}.
\end{equation}
Conversely, for any $f\in X$, we can use the above method to construct $\phi\in \mathcal{D}$ such that
\begin{equation}
\Delta \phi+\frac{\phi}{2}+\frac{x}{2}\cdot\nabla\phi-\lambda\phi=f.
\end{equation}
Thus we conclude $\mathcal{L}:\mathcal{D}\rightarrow X$ is a densely defined, closed operator with resolvent set
\begin{equation*}
\rho(\mathcal{L})\supseteq \{\lambda: {\rm Re}\lambda >-\frac{1}{4} \}.
\end{equation*}

\smallskip
\noindent
{\bf Remarks.} The above method using the relation between $\mathcal{L}$ and heat equation also provides a simple proof that Schwartz class divergence free vector fields are dense in $\mathcal{D}$, if we equip
$\mathcal{D}$ with the natural norm
\begin{equation}
\|\phi\|_{\mathcal{D}}:=\|\mathcal{L}\phi\|_X,\quad \phi\in \mathcal{D}.
\end{equation}
To prove the claim, take any $\phi\in\mathcal{D}$ and let $f=\mathcal{L}\phi\in X$. It is clear that we can take $f_n\in \mathcal{S}(R^3)$ with ${\rm div}\,f_n=0$ such that $f_n\to f$ in $X$ as $n\to \infty$. Let $\phi_n\in\mathcal{D}$ with $\mathcal{L}\phi_n=f_n$. From the heat equation formulation, it is not hard to verify that $\phi_n\in S(R^3)$. Since $\phi_n\to\phi$ in $\mathcal{D}$ as $n\to \infty$, the claim is proved.\\

We next show that $\mathcal{K}(a)$ is compact relative to $\mathcal{L}$.
\begin{lemma}\label{lm:relativecompactnessofperturbation}
Let $\mathcal{K}(a)$ and $\mathcal{L}$ be defined as above. Then $\mathcal{K}(a)$ is compact relative to $\mathcal{L}$ in the following sense. Let $\phi_n\in \mathcal{D}$ with
\begin{equation}
\|\mathcal{L}\phi_n\|_X\leq M, {\rm ~for~some~}M>0.
\end{equation}
Then the sequence $\mathcal{K}(a)\phi_n$ is compact in $X$.
\end{lemma}

\smallskip
\noindent
{\bf Proof.} By the estimates in equation (\ref{eq:estimateforlinearoperator}), we see $\phi_n,~\nabla\phi_n$ and $\nabla^2\phi_n$ are uniformly bounded in $X$. Note also that $a$ decays at spatial infinity.
From these, it's standard to show compactness.\\

Now we turn to the semigroup generated by $\mathcal{L}$ and its compact perturbations.
We will consider the semigroup generated by $\mathcal{L}-\mathcal{K}(a)$. Instead of using the various generation theorems for semigroups, we take the approach of directly solving the following initial value
problem to get a solution operator $S_a(t)$, and show the generator is $\mathcal{L}-\mathcal{K}(a)$:
\begin{eqnarray}\label{eq:semigroup}
\left\{\begin{array}{rl}
\partial_t\phi&=\Delta\phi+\frac{x}{2}\cdot\nabla\phi+\frac{\phi}{2}-a\cdot\nabla\phi-\phi\cdot\nabla a+\nabla p\\
{\rm div~}\phi&=0\\
\phi(\cdot,0)&=\phi_0
      \end{array}\right.{\rm ~in~}R^3\times(0,\infty),
\end{eqnarray}
where $\phi_0\in X$.\\

\begin{lemma}\label{lm:semigroup}
For any $\phi_0\in X$, there exists a unique function $\phi(x,t)\in L^{\infty}(0,T;X)$, with
\begin{equation}
\sup_{t\in (0,T)}\left(t\|\nabla^2\phi(\cdot,t)\|_X+t\|\partial_t\phi-\frac{x}{2}\cdot\nabla\phi(\cdot,t)\|_X\right)<\infty, {\rm ~for~any~} T>0,
\end{equation}
such that $\phi$ satisfies equation (\ref{eq:semigroup}) in the sense of distributions with appropriate distribution $p$ in $R^3\times(0,\infty)$, and
\begin{equation}
\lim_{t\to 0+}\|\phi(\cdot,t)-\phi_0\|_X=0.
\end{equation}
Moreover, we have the following estimate
\begin{equation}
\sup_{t\in[0,T)}\left(\|\phi(\cdot,t)\|_X+t\|\nabla^2\phi(\cdot,t)\|_X+t\|\partial_t\phi-\frac{x}{2}\cdot\nabla\phi(\cdot,t)\|_X\right)\leq C(\|a\|_Y,T)\|\phi_0\|, {\rm ~for~} T>0.
\end{equation}
In addition, the solution semigroup is generated by $\mathcal{L}-\mathcal{K}(a)$.
\end{lemma}

\smallskip
\noindent
 {\bf Proof.} We first prove existence. Again we use the relation with heat equation, and set
\begin{equation}
h(x,t)=\frac{1}{\sqrt{t+1}}\phi(\frac{x}{\sqrt{t+1}},\log{(t+1)}), {\rm ~for ~}t>0,~x\in R^3.
\end{equation}
It is easy to check $h$ satisfies
\begin{eqnarray}
\left.\begin{array}{rl}
     \partial_th-\Delta h+\tilde{a}\cdot\nabla h+h\cdot\nabla \tilde{a}+\nabla p&=0\\
                                                                             {\rm div~}h&=0\\
                                                                                    h(\cdot,0)&=\phi_0
 \end{array}\right\}{\rm ~in~} R^3\times(0,\infty),
\end{eqnarray}
where $\tilde{a}(x,t):=\frac{1}{\sqrt{t+1}}a(\frac{x}{\sqrt{t+1}})$. It is clear we can solve this equation with
\begin{equation}
\|h(\cdot,t)\|_X+\|\nabla h(\cdot,t)\|_Xt^{\frac{1}{2}}+t\|\nabla^2 h(\cdot,t)\|_X+t\|\partial_t h(\cdot,t)\|_X\leq C(T,\|a\|_Y)\|\phi_0\|_X {\rm ~for~} t\in [0,T).
\end{equation}
Moreover, $h(\cdot,t)$ is continuous in $t$ with values in $X$. Using the relation
\begin{equation}
\phi(y,s)=e^{\frac{s}{2}}h(e^{\frac{s}{2}}y,e^s-1), {\rm ~for~} y\in R^3,~s\ge 0.
\end{equation}
we obtain a solution $\phi$ with
\begin{equation}
\|\phi(\cdot,s)\|_X+\|\nabla\phi(\cdot,s)\|_Xs^{\frac{1}{2}}+s\|\partial_s\phi-\frac{x}{2}\cdot\nabla\phi(\cdot,s)\|_X+s\|\nabla^2\phi(\cdot,s)\|_X\leq C(T,\|a\|_Y)\|\phi_0\|_X,
\end{equation}
for $s\in[0,T)$. The uniqueness follows from the corresponding result for Stokes system.\\
Denote the solution semigroup as
\begin{equation}
S_a(s)\phi_0:=\phi(\cdot,s).
\end{equation}
Suppose
\begin{equation}
S_a(s)=e^{As},
\end{equation}
i.e, $S_a(s)$ is generated by closed operator $A$ with domain $\mathcal{D}(A)$. It is easy to see $\mathcal{D}(A)$ contains the set of Schwartz class
divergence free vector fields, and $A\phi=(\mathcal{L}-\mathcal{K}(a))\phi$ for such $\phi$. To show $A\equiv \mathcal{L}-\mathcal{K}(a)$, it suffices to show $\mathcal{D}\subseteq\mathcal{D}(A)$.
Take $\lambda>0$ sufficiently large, then $\lambda\in \rho(A)$.
Take any $\phi\in\mathcal{D}$, and take a sequence of Schwartz class divergence free $\phi_n$ such that
\begin{equation}
\|\mathcal{L}(\phi_n-\phi)\|_X\to 0, {\rm ~as~}n\to \infty.
\end{equation}
This implies $\phi_n\to\phi$ as $n\to \infty$ in $X$ by the invertibility of $\mathcal{L}$. Thus $(A-\lambda)\phi_n$ converges in $X$. Since $\phi_n$ converges in $X$, by the closedness of $A$, we see $\phi\in \mathcal{D(A)}$ and the proof is complete.\\

We have the following estimate for the semigroup.
\begin{equation}
\|e^{(\mathcal{L}-\mathcal{K}(a))t}\|+t^{\frac{1}{2}}\|\nabla e^{(\mathcal{L}-\mathcal{K}(a))t}\|+t\|\nabla^2e^{(\mathcal{L}-\mathcal{K}(a))t}\|+t\|(\partial_t-\frac{x}{2}\cdot\nabla)e^{(\mathcal{L}-\mathcal{K}(a))t}\|\leq C(T,\|a\|_Y),
\end{equation}
for $t\in[0,T)$, where the norm is operator norm from $X$ to $X$. Of particular importance is the case $a\equiv 0$. Then for $t\ge 1$,
\begin{equation}
\|h(\cdot,t)\|_{L^2(R^3)}+t^{\frac{3}{8}}\|h(\cdot,t)\|_{L^4(R^3)}\leq C\|\phi_0\|_{L^2(R^3)}\leq C\|\phi_0\|_X.
\end{equation}
Thus
\begin{equation}
\|\phi(\cdot,s)\|_X\leq Ce^{-\frac{s}{4}}\|\phi_0\|_X.
\end{equation}
In summary, we conclude the semigroup $e^{\mathcal{L}s}$ is exponentially decreasing with exponent $\frac{1}{4}$. \\

For later applications, it is essential to know the long time behavior of the semigroup $e^{(\mathcal{L}-\mathcal{K}(a))t}$. To study such behavior we need to recall some results from the general theory of semigroups from Page 248-Page 259 of \cite{Nagel}.  We begin by the following definition.
\begin{definition}
Let $A$ be a densely defined, closed operator on $X$, and denote $\sigma(A)$ the spectrum of $A$. Suppose $A$ generates a strongly continuous semigroup $e^{At}$.  Define
\begin{equation}
s(A):=\sup\{{\rm Re}\,\lambda: \lambda\in \sigma(A)\};
\end{equation}
define the growth bound for $e^{At}$:
\begin{equation}
\omega_0:=\inf\left\{\omega\in R: {\rm ~there~exists~}M_{\omega}\ge 1 {\rm ~such~that~} \|e^{At}\|\leq M_{\omega}e^{\omega t} {\rm ~for~all~}t\ge0\right\};
\end{equation}
for a bounded linear operator $T:X\to X$, define
\begin{equation}
\|T\|_{ess}:=\inf\left\{\|T-K\|: K {\rm ~is~compact~from~}X~{\rm to}~X \right\};
\end{equation}
define the essential growth bound for $e^{At}$ as
\begin{equation}
\omega_{ess}:=\inf_{t>0}\frac{1}{t}\log{\|e^{At}\|_{ess}}.
\end{equation}
\end{definition}
Let
\begin{equation}
\sigma_{ess}(A):=\{\lambda\in\mathbb{C}: A-\lambda I {\rm ~is~not~Fredholm}\},
\end{equation}
and
\begin{equation}
r_{ess}(T):=\sup\{|\lambda|:\lambda\in \sigma_{ess}(T)\}.
\end{equation}
We recall the following lemmas from \cite{Nagel} (Proposition 2.2 in Page 251, Proposition 2.10 and Corollary 2.11 in page 258).
\begin{lemma}\label{lm:formulaforgrowthbound}
\begin{equation}
\omega_0=\inf_{t>0}\frac{1}{t}\log{\|e^{At}\|}=\lim_{t\to\infty}\frac{1}{t}\log{\|e^{At}\|}=\frac{1}{t_0}\log{r(e^{At_0})},
\end{equation}
for each $t_0>0$, where $r(T)$ is the spectral radius of $T$.
\end{lemma}

\begin{lemma}\label{lm:essentialgrowthbound1}
\begin{equation}
\omega_{ess}=\lim_{t\to \infty}\frac{1}{t}\log{\|e^{At}\|_{ess}}=\frac{1}{t_0}\log{r_{ess}(e^{At_0})}\leq \omega_0<\infty,
\end{equation}
for each $t_0>0$.
\end{lemma}

\begin{lemma}\label{lm:generalsemigroupstuff}
Let $e^{At}$ be a strongly continous semigroup on the Banach space $X$ with generator $A$. Then
\begin{equation}
\omega_0=\max\{\omega_{ess}, ~s(A)\}.
\end{equation}
Moreover, for every $\omega>\omega_{ess}$ the set $\sigma_c:=\sigma(A)\cap\{\lambda\in\mathbb{C}:{\rm ~Re\,\lambda\ge\omega}\}$ is finite and the
corresponding spectral projection has finite rank.
\end{lemma}
Note that $\omega_{ess}$ is stable under compact perturbations, and $s(A)$ is usually easy to understand, thus the above lemma provides us with a powerful tool to control the growth bound.
Now we are ready to estimate the growth bound for semigroup $e^{(\mathcal{L}-\mathcal{K}(a))t}$.
\begin{lemma}\label{lm:essentialgrowthbound2}
Let $\mathcal{L}$ and $\mathcal{K}(a)$ be defined as above. Then
\begin{equation}
\omega_{ess}(e^{(\mathcal{L}-\mathcal{K}(a))t})=\omega_{ess}(e^{\mathcal{L}t})\leq \omega_0(e^{\mathcal{L}t})\leq -\frac{1}{4}.
\end{equation}
Thus
\begin{equation}
\omega_0(e^{(\mathcal{L}-\mathcal{K}(a))t})\leq \max\left\{-\frac{1}{4},~s(\mathcal{L}-\mathcal{K}(a))\right\}.
\end{equation}
Moreover, for any $\omega>-\frac{1}{4}$, the set $\{\lambda\in \sigma(\mathcal{L}-\mathcal{K}(a)): {\rm ~Re}\,\lambda\ge\omega\}$ is finite.
\end{lemma}

\smallskip
\noindent
{\bf Proof.} By the results in Lemma \ref{lm:essentialgrowthbound1} to show
\begin{equation*}
\omega_{ess}(e^{(\mathcal{L}-\mathcal{K}(a))t})=\omega_{ess}(e^{\mathcal{L}t}),
\end{equation*}
 we only need to show
\begin{equation*}
r_{ess}(e^{(\mathcal{L}-\mathcal{K}(a))t})=r_{ess}(e^{\mathcal{L}t}).
\end{equation*}
Recall the following formula relating $e^{(\mathcal{L}+\mathcal{K}(a))t}$ and $e^{\mathcal{L}t}$:
\begin{equation}\label{eq:formularelatingtwosemigroups}
e^{\mathcal{L}-\mathcal{K}(a)}\phi-e^{\mathcal{L}}\phi=-\int_0^1e^{\mathcal{L}(1-s)}\mathcal{K}(a)e^{(\mathcal{L}-\mathcal{K}(a))s}\phi ds, {\rm ~for~}\phi\in X.
\end{equation}
The above formula can be obtained by differentiating $e^{\mathcal{L}(1-s)}e^{(\mathcal{L}-\mathcal{K}(a))s}\phi$ in $s$ and integrate from $0$ to $1$, in the case $\phi\in \mathcal{D}$. The general case then follows from approximation.
Since the essential spectral radius is invariant with respect to compact perturbation, we only need to verify the right hand side of equation (\ref{eq:formularelatingtwosemigroups}) is a compact operator from $X$ to $X$. Take $\phi_n\in X$ with $\|\phi_n\|_X\leq M$ and $\phi_n\rightharpoonup \phi\in X$. By the estimates of the semigroup and the decay estimates, we see for $s>0$,
\begin{equation*}
\|\mathcal{K}(a)e^{(\mathcal{L}-\mathcal{K}(a))s}(\phi_n-\phi)\|_X\to 0.
\end{equation*}
Moreover, we have
\begin{equation}
\|\mathcal{K}(a)e^{(\mathcal{L}-\mathcal{K}(a))s}(\phi_n-\phi)\|_X\leq C(\|a\|_Y)Ms^{-\frac{1}{2}}.
\end{equation}
Thus
\begin{eqnarray*}
&&\|\int_0^1e^{\mathcal{L}(1-s)}\mathcal{K}(a)e^{(\mathcal{L}-\mathcal{K}(a))s}(\phi_n-\phi) ds\|_X\\
&&\quad\quad\leq C\int_0^1\|\mathcal{K}(a)e^{(\mathcal{L}-\mathcal{K}(a))s}(\phi_n-\phi)\|_Xds\to 0, {\rm ~as~}n\to\infty,
\end{eqnarray*}
by dominated convergence theorem. Therefore the compactness is proved. The other claims in the lemma now follows directly from Lemma \ref{lm:generalsemigroupstuff}.\\

We shall need the following simple continuity result for the semigroup $e^{(\mathcal{L}-\mathcal{K}(a))t}$.
\begin{lemma}\label{lm:parametercontinuityofsemigroup}
For any $M>0$, let $B_M(0)\subset Y$ be the ball with radius $M$ centered at origin in $Y$. Then for any $t>0$ and $\epsilon>0$, there exists $\delta=\delta(M,t,\epsilon)>0$ such that if $a,\,b\in B_M(0)$ and
\begin{equation}
\|a-b\|_{C^1(R^3)}<\delta,
\end{equation}
then
\begin{equation}
\|e^{(\mathcal{L}-\mathcal{K}(a))t}-e^{(\mathcal{L}-\mathcal{K}(b))t}\|<\epsilon,
\end{equation}
where the norm is operator norm from $X$ to $X$.
\end{lemma}

\smallskip
\noindent
{\bf Proof.} We have the following formula relating the semigroups $e^{(\mathcal{L}-\mathcal{K}(a))t}$ and $e^{(\mathcal{L}-\mathcal{K}(b))t}$:
\begin{equation}
e^{(\mathcal{L}-\mathcal{K}(a))t}-e^{(\mathcal{L}-\mathcal{K}(b))t}=\int_0^te^{(\mathcal{L}-\mathcal{K}(b))(t-s)}\mathcal{K}(b-a)e^{(\mathcal{L}-\mathcal{K}(a))s}ds.
\end{equation}
Thus for any $\phi\in X$:
\begin{eqnarray*}
&&\|e^{(\mathcal{L}-\mathcal{K}(a))t}\phi-e^{(\mathcal{L}-\mathcal{K}(b))t}\phi\|_X\\
&&\quad\quad \leq C(M,t)\int_0^t\|a-b\|_{C^1(R^3)}\max\{s^{-\frac{1}{2}},1\}\|\phi\|_Xds\\
&&\quad\quad\quad\leq C(M,t)\|a-b\|_{C^1(R^3)}\|\phi\|_X.
\end{eqnarray*}
The lemma follows from the above inequality easily.\\

Lemma \ref{lm:essentialgrowthbound2} gives an estimate of the exponent of exponential growth for the semigroup $e^{(\mathcal{L}-\mathcal{K}(a))t}$. For any $\beta>\omega_0$, there exists $C>0$ such that
\begin{equation}
\|e^{(\mathcal{L}-\mathcal{K}(a))t}\|_{X\to X}\leq Ce^{\beta t}.
\end{equation}
However it is not clear what the value of $C$ depends on. In this regard, we give the following simple lemma.
\begin{lemma}\label{lm:quantitativeestimate}
For $M>0$, $\omega_{\ast}\in R$, let $a\in Y$, $\omega_0$ be defined as above with $\|a\|_Y\leq M$, $\omega_0\leq \omega_{\ast}$. Then for any $\beta>\max\{\omega_{\ast},-\frac{1}{4}\}$, there exists a constant $C=C(M,\beta,\omega_{\ast})$ such that
\begin{equation}
\|e^{(\mathcal{L}-\mathcal{K}(a))t}\|_{X\to X}\leq Ce^{\beta t}.
\end{equation}
\end{lemma}

\smallskip
\noindent
{\bf Proof.} The proof is a standard application of compactness argument and the continuous dependence of the semigroup $e^{(\mathcal{L}-\mathcal{K}(a))t}$ on $a$ in a weaker space ($C^1(R^3)$) in which $Y$ is compact. Suppose the lemma is not true. Then there exists $a_n\in Y$ with $\|a_n\|_Y\leq M$ and $\omega_0(e^{(\mathcal{L}-\mathcal{K}(a_n))t})\leq \omega_{\ast}$, such that
\begin{equation}
\sup_{t\ge 0}e^{-\beta t}\|e^{(\mathcal{L}-\mathcal{K}(a_n))t}\|_{X\to X}\ge n, {\rm ~for~each~integer~}n>0.
\end{equation}
We can pass to a subsequence (and still denoting the new sequence as $a_n$) such that $a_n\to a$ in $C^1(R^3)$, for some $a\in Y$. Take $\beta>\beta_1>\max\{\omega_{\ast},-\frac{1}{4}\}$,  by Lemma \ref{lm:essentialgrowthbound2} we know
\begin{equation}
\sigma(e^{(\mathcal{L}-\mathcal{K}(a))t})\subseteq B_{e^{\beta_1t}}\cup \{~{\rm finitely~many~eigenvalues~with~}|\lambda|>e^{\beta_1 t}\}.
\end{equation}
Since $e^{(\mathcal{L}-\mathcal{K}(a_n))t}$ converges to $e^{(\mathcal{L}-\mathcal{K}(a))t}$ as $n\to \infty$ by Lemma \ref{lm:parametercontinuityofsemigroup}, and $e^{(\mathcal{L}-\mathcal{K}(a_n))t}$ has no eigenvalues in $\{\lambda\in \mathbb{C}|~|\lambda|>e^{\beta_1t}\}$, we conclude $e^{(\mathcal{L}-\mathcal{K}(a))t}$ has no eigenvalues in $\{\lambda\in \mathbb{C}|~|\lambda|>e^{\beta_1t}\}$ as well. Thus
\begin{equation}
\sigma\left(e^{(\mathcal{L}-\mathcal{K}(a))t}\right)\subseteq B_{e^{\beta_1t}}.
\end{equation}
Thus $\omega_0(e^{(\mathcal{L}-\mathcal{K}(a))t})<\beta$. Therefore by Lemma \ref{lm:formulaforgrowthbound} we have
\begin{equation}
\frac{1}{t_0}\log{\|e^{(\mathcal{L}-\mathcal{K}(a))t_0}\|}<\beta,
\end{equation}
for some $t_0$ sufficiently large. Thus
\begin{equation}
\|e^{(\mathcal{L}-\mathcal{K}(a))t_0}\|<e^{\beta t_0}.
\end{equation}
By Lemma \ref{lm:parametercontinuityofsemigroup}, we see for $n\ge n_0$ ($n_0$ sufficiently large):
\begin{equation}
\|e^{(\mathcal{L}-\mathcal{K}(a_n))t_0}\|<e^{\beta t_0}.
\end{equation}
Denote
\begin{equation}
\gamma(M,t_0):=\sup_n\sup_{t\in(0,t_0)}\|e^{(\mathcal{L}-\mathcal{K}(a_n))t}\|<\infty,
\end{equation}
we easily obtain
\begin{equation}
\|e^{(\mathcal{L}-\mathcal{K}(a_n))t}\|\leq \gamma(M,t_0)C(t_0,\beta)e^{\beta t}, {\rm ~for~}n\ge n_0, ~t>0.
\end{equation}
This is a contradiction. The lemma is proved.
\end{section}

\begin{section}{Generalized Stokes system with singular lower order terms}
In order to localize the forward self similar solutions (which have infinite energy), we need to consider perturbation of these solutions in regular spaces. The following lemma is important in such considerations.
\begin{lemma}\label{lm:solvabilityofsingularstokes}
Let $a\in Y$ be such that $\beta:=s(\mathcal{L}-\mathcal{K}(a))<\frac{1}{8}$. Denote $\tilde{a}(x,t)=\frac{1}{\sqrt{t}}a(\frac{x}{\sqrt{t}})$. Let $f$ be such that
\begin{eqnarray}
&&M:=\esssup_{t\in(0,1)}\left(\|f(\cdot,t)\|_{L^4(R^3)}t+t^{\frac{5}{8}}\|f(\cdot,t)\|_{L^2(R^3)}\right)<\infty, {\rm ~and}\\
&&\lim_{t\to0+}\left(\|f(\cdot,t)\|_{L^4(R^3)}t+t^{\frac{5}{8}}\|f(\cdot,t)\|_{L^2(R^3)}\right)=0.
\end{eqnarray}
Then there exists a unique solution $u\in L_t^{\infty}L_x^4(R^3\times (0,1))$ to the Stokes system with singular lower order terms
\begin{eqnarray}\label{eq:singularstokesequation}
\left.\begin{array}{rl}
    \partial_tu-\Delta u+\tilde{a}\cdot\nabla u+u\cdot\nabla\tilde{a}+\nabla p&=f\\
                                                                  {\rm div~}u&=0\\
         \end{array}\right\} {\rm ~in~}R^3\times(0,1),
\end{eqnarray}
with divergence free initial data $u(\cdot,0)=u_0\in L^4(R^3)$ in the sense that
\begin{equation}
\lim_{t\to 0+}\|u(\cdot,t)-u_0\|_{L^4(R^3)}=0.
\end{equation}
Moreover, $u$ satisfies
\begin{equation}
\|u\|_{L_t^{\infty}L_x^4(R^3\times(0,1))}+\sup_{t\in(0,1)}t^{\frac{1}{2}}\|\nabla u(\cdot,t)\|_{L_x^4(R^3)}\leq C(\|a\|_Y,\beta)\left(\|u_0\|_{L_x^4(R^3)}+M\right),
\end{equation}
and
\begin{equation}
\lim_{t\to 0+}t^{\frac{1}{2}}\|\nabla u(\cdot,t)\|_{L^4(R^3)}=0.
\end{equation}
\end{lemma}

\smallskip
\noindent
{\bf Proof.} We first consider the case $u_0\in C_c^{\infty}(R^3)$ and $f\in C_c^{\infty}(R^3\times(0,1))$ with $\|u_0\|_{L^4(R^3)}+M=1$. To construct a solution, we write
\begin{equation}\label{eq:singularStokes}
u(x,t):=\phi(\frac{x}{\sqrt{t}},\log{t})+e^{\Delta t}u_0, {\rm ~for~}(x,t)\in R^3\times(0,1),
\end{equation}
and
\begin{equation}
f(x,t)=\frac{1}{t}g(\frac{x}{\sqrt{t}},\log{t}), ~(x,t)\in R^3\times (0,1).
\end{equation}
Denote $e^{\Delta t}u_0(x)=b(\frac{x}{\sqrt{t}},\log{t})$ for $(x,t)\in R^3\times(0,1)$. Since $u$ must satisfy equation (\ref{eq:singularstokesequation}), we obtain the following equation for $\phi$
\begin{eqnarray}\label{eq:derivedequation}
\left\{\begin{array}{rl}
        \partial_s\phi&=(\mathcal{L}-\mathcal{K}(a)-\frac{1}{2})\phi-(a\cdot\nabla b+b\cdot\nabla a)+g+\nabla p\\
          {\rm div~}\phi&=0\\
       \end{array}\right.{\rm ~for~}(y,s)\in R^3\times(-\infty,0).
\end{eqnarray}
We can explicitly write down a solution as
\begin{equation}\label{eq:solutionformula}
\phi(\cdot,s):=-\int_{-\infty}^se^{(\mathcal{L}-\mathcal{K}(a)-\frac{1}{2})(s-\tau)}\mathbb{P}(a\cdot\nabla b+b\cdot\nabla a-g)(\cdot,\tau)d\tau,
\end{equation}
where $\mathbb{P}$ is Helmholtz projection to divergence free vector fields. To justify the formula, we must obtain appropriate estimates. Firstly, note
that (denoting $h(x,t)=e^{\Delta t}u_0(x)$)
\begin{equation}
b(y,s)=h(e^{\frac{s}{2}}y,e^s).
\end{equation}
Thus by scale-invariance of heat equation, we see $b(y,s)$ is the solution to heat equation at time $1$ with initial data $u_0(e^{\frac{s}{2}}\cdot)$. Therefore we obtain the following estimates for $b$
\begin{equation}
\|b(\cdot,s)\|_{L_y^4\cap L_y^{\infty}(R^3)}+\|\nabla b(\cdot,s)\|_{L_y^4\cap L_y^{\infty}(R^3)}\leq C e^{-\frac{3s}{8}}.
\end{equation}
Thus from the assumptions on $a$, we get
\begin{equation}
\|\mathbb{P}(a\cdot\nabla b+b\cdot\nabla a)(\cdot,s)\|_X\leq C(\|a\|_Y) e^{-\frac{3s}{8}}.
\end{equation}
Note
\begin{equation}
g(y,s)=e^sf(e^{\frac{s}{2}}y,e^s), {\rm ~for~}(y,s)\in R^3\times(-\infty,0),
\end{equation}
thus by the assumption on $f$, we obtain
\begin{equation}
\|g(\cdot,s)\|_X\leq Ce^{-\frac{3s}{8}}, \, {\rm a\,e}\,s\in (-\infty,0).
\end{equation}
By Lemma \ref{lm:essentialgrowthbound2}
\begin{equation}
\omega_0(e^{(\mathcal{L}-\mathcal{K}(a))t})\leq\max\{\beta,-\frac{1}{4}\}<\frac{1}{8}.
\end{equation}
Take $\max\{\beta,-\frac{1}{4}\}<\beta_1<\frac{1}{8}$, by the definition of $\omega_0$ and Lemma \ref{lm:quantitativeestimate}, there exists constant $C(\|a\|_Y,\beta)>0$ such that
\begin{equation}\label{eq:decayintime}
\|e^{(\mathcal{L}-\mathcal{K}(a))t)}\|_{X\to X}\leq C(\|a\|_Y,\beta) e^{\beta_1t} {\rm ~for~any~}t\ge 0.
\end{equation}
Combining this decay estimate with the local regularity property of the semigroup, we see also
\begin{equation}
\|\nabla e^{(\mathcal{L}-\mathcal{K}(a))t}\|_{X\to X}\leq C(\|a\|_Y,\beta)\max\{\frac{1}{\sqrt{t}},1\}e^{\beta_1t}.
\end{equation}
Thus we get
\begin{equation}
\|\phi(\cdot,s)\|_X\leq \int_{-\infty}^sC(\|a\|_Y,\beta)e^{(\beta_1-\frac{1}{2})(s-\tau)}e^{-\frac{3\tau}{8}}d\tau\leq C(\|a\|_Y,\beta)e^{-\frac{3s}{8}},
\end{equation}
where we used the condition $\beta_1<\frac{1}{8}$. Similarly
\begin{equation}
\|\nabla\phi(\cdot,s)\|_X\leq C(\|a\|_Y,\beta)e^{-\frac{3s}{8}}.
\end{equation}
Consequently, we have
\begin{equation}
\sup_{t\in(0,1)}\|u(\cdot,t)-e^{\Delta t}u_0\|_{L^4(R)}\leq t^{\frac{3}{8}}\|\phi(\cdot,\log{t})\|_{L_x^4(R^3)}\leq C(\|a\|_Y,\beta),
\end{equation}
and
\begin{equation}
\sup_{t\in(0,1)}t^{\frac{1}{2}}\|\nabla u(\cdot,t)\|_{L^4(R^3)}\leq C(\|a\|_Y,\beta).
\end{equation}
Note that these estimates do not depend on the assumption that $u_0\in C_c^{\infty}(R^3)$ and $f\in C_c^{\infty}(R^3\times(0,1))$. Now with this condition, one can check
\begin{equation}
\|\mathbb{P}(a\cdot\nabla b+b\cdot\nabla a)(\cdot,s)\|_X+\|g(\cdot,s)\|_X\leq C(u_0)e^{-\delta s} {\rm ~for~some~}\delta\in (0,\frac{3}{8}),
\end{equation}
and consequently
\begin{equation}
\|\phi(\cdot,s)\|_X+\|\nabla\phi(\cdot,s)\|_X\leq C(\|a\|_Y,u_0,\beta)e^{-\delta s}.
\end{equation}
In this case, we can show by direct differentiation that $\phi$ satisfies equation (\ref{eq:derivedequation}), and using the relation $u(x,t)=\phi(\frac{x}{\sqrt{t}},\log{t})+e^{\Delta t}u_0$, it is straightforward to verify that $u$ satisfies equation (\ref{eq:singularstokes}) with
\begin{eqnarray*}
&&\|u(\cdot,t)-e^{\Delta t}u_0\|_{L^4(R)}+t^{\frac{1}{2}}\|\nabla u(\cdot,t)\|_{L^4(R^3)}\\
&&\quad\quad\leq t^{\frac{3}{8}}\left(\|\phi(\cdot,\log{t})\|_{L_x^4(R^3)}+\|\nabla \phi(\cdot,\log{t})\|_{L^4(R^3)}\right)+t^{\frac{1}{2}}\|\nabla e^{\Delta t}u_0\|_{L^4(R^3)}\\
&&\quad\quad\quad\leq t^{\frac{3}{8}-\delta}C(u_0,\|a\|_Y,\beta)+t^{\frac{1}{2}}\|\nabla e^{\Delta t}u_0\|_{L^4(R^3)}\to 0,
\end{eqnarray*}
as $t\to 0+$. This completes the proof when the initial data $u_0\in C_c^{\infty}$, the general case then follows from approximation. \\
To prove uniqueness, suppose $u$ is a distributional solution with $u_0\equiv 0$ and $f\equiv 0$, then
\begin{equation}
u(\cdot,t)=-\int_0^te^{\Delta(t-s)}\mathbb{P}{\rm div}\,(\tilde{a}\otimes u+u\otimes\tilde{a})(\cdot,s)ds.
\end{equation}
From this we can obtain
\begin{equation}
\|u(\cdot,t)\|_{L^2(R^3)}\leq Ct^{\frac{3}{8}}.
\end{equation}
Write $u(x,t)=\phi(\frac{x}{\sqrt{t}},\log{t})$. Since for $t$ positive $u(x,t)$ is regular,  $\phi$ is given by the semigroup
\begin{equation}
\phi(\cdot,s)=e^{(\mathcal{L}-\mathcal{K}(a))(s+M)}\phi(\cdot,-M).
\end{equation}
Note
\begin{equation}
\|\phi(\cdot,-M)\|_X\leq Ce^{\frac{3}{8}M},
\end{equation}
by the estimates on $u$ and the relation between $u$ and $\phi$. Thus by the decay property of the semigroup, we get
\begin{equation}
\|\phi(\cdot,s)\|_X\leq C(\|a\|_Y,\beta)e^{(\beta_1-\frac{1}{2})(s+M)}e^{\frac{3}{8}M}\leq C(\|a\|_Y,\beta,s)e^{(\beta_1-\frac{1}{8})M}\to 0, {\rm ~as~}M\to \infty+.
\end{equation}
Thus $\phi\equiv 0$ and consequently $u\equiv 0$.\\

Now we are ready to prove the first main theorem of this paper.
\begin{theorem}\label{th:perturbingsingularsolution}
Let $a\in Y$ be such that $\frac{1}{8}>\beta:=s(\mathcal{L}-\mathcal{K}(a))$ which is the maximum of real part of any eigenvalue of $\mathcal{L}-\mathcal{K}(a)$. Denote $\tilde{a}(x,t)=\frac{1}{\sqrt{t}}a(\frac{x}{\sqrt{t}})$. Let
\begin{equation}
\tilde{b}(x,t):=\frac{1}{\sqrt{t}}b(\frac{x}{\sqrt{t}},\log{t}), {\rm ~for~}x\in R^3,~t\in(0,1).
\end{equation}
Suppose
\begin{equation}\label{eq:assumptionofb}
\epsilon=\esssup_{s\in(-\infty,0)}(\|b(\cdot,s)\|_X+\|\nabla b(\cdot,s)\|_X)
\end{equation}
is sufficiently small depending on $\|a\|_Y,\,\beta$. Let $T=T(\|a\|_Y,\beta,\|u_0\|_{L^4(R^3)})>0$ be sufficiently small.
Then there exists a unique solution $u\in Z_T$ to the generalized Navier Stokes system with singular lower order terms
\begin{eqnarray}\label{eq:singularstokes}
\left.\begin{array}{rl}
    \partial_tu-\Delta u+\tilde{a}\cdot\nabla u+u\cdot\nabla\tilde{a}+\tilde{b}\cdot\nabla u+u\cdot\nabla \tilde{b}+u\cdot\nabla u+\nabla p&=0\\
                                                                  {\rm div~}u&=0\\
         \end{array}\right\} {\rm ~in~}R^3\times(0,T),
\end{eqnarray}
with initial data $u(\cdot,0)=u_0\in L^4(R^3)$ in the sense that
\begin{equation}
\lim_{t\to0+}\|u(\cdot,t)-u_0\|_{L^4(R^3)}=0.
\end{equation}
In the above,
\begin{equation}
Z_T:=\{u\in L_t^{\infty}L_x^4(R^3\times(0,T)): ~\esssup_{t\in (0,T)}t^{\frac{1}{2}}\|\nabla u(\cdot,t)\|_{L^4(R^3)}<\infty, ~\lim_{t\to0+}t^{\frac{1}{2}}\|\nabla u(\cdot,t)\|_{L^4(R^3)}=0\},
\end{equation}
equipped with the natural norm
\begin{equation}
\|u\|_{Z_T}:=\esssup_{t\in (0,T)}\left(\|u(\cdot,t)\|_{L^4(R^3)}+t^{\frac{1}{2}}\|\nabla u(\cdot,t)\|_{L^4(R^3)}\right).
\end{equation}
Moreover, $u$ satisfies
\begin{equation}
\|u\|_{L_t^{\infty}L_x^4(R^3\times(0,T))}+\sup_{t\in(0,T)}t^{\frac{1}{2}}\|\nabla u(\cdot,t)\|_{L_x^4(R^3)}\leq C(\beta,\|a\|_Y)\|u_0\|_{L_x^4(R^3)}.
\end{equation}
\end{theorem}

\smallskip
\noindent
{\bf Proof.} Denote $T(t)u_0+B(f)(t)$ as the solution to equation (\ref{eq:singularstokesequation}) (we suppress the dependence on $a$).
Then by Lemma \ref{lm:solvabilityofsingularstokes}, we have
\begin{eqnarray}
&&\|T(t)u_0\|_{Z_T}\leq C(\|a\|_Y,\beta)\|u_0\|_{L^4(R^3)},\\
&&\|B(f)(t)\|_{Z_T}\leq C(\|a\|_Y,\beta)\esssup_{t\in(0,T)}\left(\|f(\cdot,t)\|_{L^4(R^3)}t+t^{\frac{5}{8}}\|f(\cdot,t)\|_{L^2(R^3)}\right).
\end{eqnarray}
To prove the theorem, we write
\begin{equation}\label{eq:contraction}
u(\cdot,t)=T(t)u_0-B(\tilde{b}\cdot\nabla u+u\cdot\nabla \tilde{b}+u\cdot\nabla u)(t),
\end{equation}
and use the contraction mapping theorem in $Z_T$ to prove existence and uniqueness once $T$ is sufficiently small. The method is standard, and we only sketch some of the details below.
By the definition of $\tilde{b}$ and the estimate \ref{eq:assumptionofb} of $b$, we obtain for $t\in(0,1)$
\begin{eqnarray}
\|\tilde{b}(\cdot,t)\|_{L^4(R^3)}&\leq&t^{-\frac{1}{8}}\epsilon,\\
\|\nabla\tilde{b}(\cdot,t)\|_{L^4(R^3)}&\leq&t^{-\frac{5}{8}}\epsilon.
\end{eqnarray}
By the following interpolation inequality:
\begin{equation}
\|f\|_{L^{\infty}(R^3)}\leq C\|f\|_{L^4(R^3)}^{\frac{1}{4}}\|\nabla f\|_{L^4(R^3)}^{\frac{3}{4}},
\end{equation}
we obtain
\begin{equation}
\|\tilde{b}(\cdot,t)\|_{L^{\infty}(R^3)}\leq Ct^{-\frac{1}{2}}\epsilon.
\end{equation}
Thus for $t\leq T$,
\begin{eqnarray*}
t\|\tilde{b}\cdot\nabla u(\cdot,t)\|_{L^4(R^3)}&\leq&Ct\|\tilde{b}(\cdot,t)\|_{L^{\infty}(R^3)}\|\nabla u(\cdot,t)\|_{L^4(R^3)}\leq C\epsilon \|u\|_{Z_T},\\
t^{\frac{5}{8}}\|\tilde{b}\cdot\nabla u(\cdot,t)\|_{L^2(R^3)}&\leq&Ct^{\frac{5}{8}}\|\tilde{b}(\cdot,t)\|_{L^{4}(R^3)}\|\nabla u(\cdot,t)\|_{L^4(R^3)}\leq C\epsilon \|u\|_{Z_T}.
\end{eqnarray*}
By similar calculations and Lemma \ref{lm:solvabilityofsingularstokes}, we can verify:
\begin{eqnarray*}
&&\|B(\tilde{b}\cdot \nabla u)\|_{Z_T}\leq C(\|a\|_Y,\beta)\epsilon \|u\|_{Z_T},\\
&&\|B(u\cdot\nabla\tilde{b})\|_{Z_T}\leq C(\|a\|_Y,\beta)\epsilon \|u\|_{Z_T},\\
&&\|B(u\cdot\nabla v)\|_{Z_T}+\|B(v\cdot\nabla u)\|_{Z_T} \leq C(\|a\|_Y,\beta)T^{\frac{1}{8}}\|u\|_{Z_T}\|v\|_{Z_T}.
\end{eqnarray*}
Denote $M:=\|T(t)u_0\|_{Z_1}$. Now take $\epsilon=\epsilon(\|a\|_Y,\beta)>0$ and $T>0$ sufficiently small, such that
\begin{eqnarray*}
&&C(\|a\|_Y,\beta)\epsilon<\frac{1}{16},\\
&&C(\|a\|_Y,\beta)MT^{\frac{1}{8}}<\frac{1}{8}.
\end{eqnarray*}
Under these conditions, it is simple to verify the right hand side of equation (\ref{eq:contraction}) is a contraction mapping from $B_{2M}\subset Z_T$ to itself, thus the theorem is proved.
\end{section}

\begin{section}{Proof of bifurcations}
Firstly, under the spectral condition ({\bf A}) we prove that we can continue the solution curve $U_{\sigma}$ as a regular curve of $\sigma$. The possibility of the continuation of the solution curve $U_{\sigma}$ is a simple application of Implicit function theorem or a direct contraction mapping argument. More precisely we have the following result.
\begin{lemma}\label{lm:continuabilityofsolutioncurve}
Suppose the spectral assumption ({\bf A}) holds. Then for $\sigma_1>\sigma_0$ sufficiently close to $\sigma_0$ and $\sigma\in(\sigma_0,\sigma_1)$, there exists $U_{\sigma}\in C^{\infty}(R^3)$ satisfying
\begin{equation}\label{eq:forwardselfsimilarequation1}
\Delta U_{\sigma}+\frac{x}{2}\cdot\nabla U_{\sigma}+\frac{1}{2}U_{\sigma}-U_{\sigma}\cdot\nabla U_{\sigma}+\nabla P=0, {\rm ~in~}R^3
\end{equation}
and
\begin{equation}\label{eq:forwardselfsimilarequation2}
\|\nabla^{\alpha}(U_{\sigma}-\sigma e^{\Delta }u_0)\|\leq \frac{C(\alpha,\sigma,u_0)}{(1+|x|)^{3+|\alpha|}}, {\rm ~for~any~multi-index~}\alpha.
\end{equation}
\end{lemma}

\smallskip
\noindent
{\bf Remark.} Since we fix $u_0$ from the beginning, in the texts below we omit the dependence of various constants on $u_0$.

\smallskip
\noindent
{\bf Proof.} Denote $U=e^{\Delta}u_0$. To prove the theorem, we seek solution $U_{\sigma}$ in the form
\begin{equation}\label{eq:aformula}
U_{\sigma}=U_{\sigma_0}+(\sigma-\sigma_0)U+\phi, {\rm ~with~}\phi\in X.
\end{equation}
Then $\phi$ satisfies $F(\phi,\sigma)=0$, where
\begin{eqnarray*}
&F(\phi,\sigma):=\Delta \phi+\frac{x}{2}\cdot\nabla \phi+\frac{\phi}{2}-U_{\sigma_0}\cdot\nabla\phi-\phi\cdot\nabla U_{\sigma_0}-(\sigma-\sigma_0)U\cdot \nabla\phi-(\sigma-\sigma_0)\phi\cdot\nabla U\\
&\quad\quad -(\sigma-\sigma_0)U\cdot\nabla U_{\sigma_0}-(\sigma-\sigma_0)U_{\sigma_0}\cdot\nabla U-(\sigma-\sigma_0)^2U\cdot\nabla U-\phi\cdot\nabla\phi+\nabla P.
\end{eqnarray*}
It is easy to verify that $F:\,\mathcal{D}\times (\sigma_0-\epsilon,\sigma_0+\epsilon)\rightarrow X$ is smooth, where $\epsilon>0$ is a small number. Note
\begin{equation}
\partial_{\phi}F(0,\sigma_0)\phi=\mathcal{L}_{\sigma_0}\phi, {\rm ~for~any~}\phi\in\mathcal{D}.
\end{equation}
By assumption ({\bf A}) $\mathcal{L}_{\sigma_0}$ is invertible. Thus,
\begin{eqnarray}
&&F(0,\sigma_0)=0\\
&&\partial_{\phi}F(0,\sigma_0) {\rm ~is~invertible~}.
\end{eqnarray}
Therefore we can apply Implicit function theorem (see for example Theorem I.1.1 in \cite{Hans}) and conclude that there exists unique $\phi(\cdot,\sigma)$ continuously differentiable with respect to $\sigma$ (with value in $\mathcal{D}$) in a small neighborhood of $(0,\sigma_0)$. Then we can verify that with $U_{\sigma}$ given by equation (\ref{eq:aformula}), $\frac{1}{\sqrt{t}}U_{\sigma}(\frac{x}{\sqrt{t}})$ is a forward self similar solution to Navier Stokes equation with initial data $\sigma u_0$. The estimates (\ref{eq:forwardselfsimilarequation2}) follow from Theorem 4.1 of \cite{JiaSverak}. The Lemma is proved.\\

Now let us consider the operator
\begin{equation}
\mathcal{L}_{\sigma}\phi=\Delta \phi+\frac{x}{2}\cdot\nabla\phi+\frac{\phi}{2}-U_{\sigma}\cdot\nabla\phi-\phi\cdot\nabla U_{\sigma}+\nabla P
\end{equation}
as defined in the introduction for $\sigma\in (\sigma_0,\sigma_1)$. By Proposition I.7.2 in \cite{Hans}, the eigenvalue $\lambda_1(\sigma)$ is differentiable in $\sigma$. By the spectral assumption ({\bf A}), for $\sigma>\sigma_0$ and sufficiently close to $\sigma_0$, ${\rm Re}\,\lambda_1(\sigma)\in (0,\frac{1}{32})$. We construct the following ancient solution.
\begin{theorem}
Suppose the spectral assumption ({\bf A}) holds. Let $\mathcal{L}_{\sigma}$ be defined as above. Then there exists $\sigma_1>\sigma_0$ sufficiently close to $\sigma_0$ such that for $\sigma\in(\sigma_0,\sigma_1)$ we have
\begin{equation}
\sigma(\mathcal{L}_{\sigma})\subset \{\lambda\in\mathbb{C}:\,{\rm Re}\,\lambda<-\frac{3}{4}\delta\}\cup \{\lambda_1(\sigma),\,\overline{\lambda_1(\sigma)},{\rm~with~Re}\,\lambda_1(\sigma)\in (0,\frac{1}{32})\}.
\end{equation}
Denoting $\beta={\rm Re}\,\lambda_1(\sigma)$, then for any $\epsilon>0$ we can find $\phi=\phi(\cdot,t,\sigma)$ with $\sigma\in (\sigma_0,\sigma_1)$ such that
\begin{equation}\label{eq:ancientequation}
\left\{\begin{array}{rl}
\partial_t\phi&=\mathcal{L}_{\sigma}\phi-\phi\cdot\nabla\phi+\nabla p\\
{\rm div~}\phi&=0
\end{array}\right.{\rm ~for~}(x,t)\in R^3\times(-\infty,0),
\end{equation}
with
\begin{equation}\label{eq:smallancientsolution}
\esssup_{t\in (-\infty,0)}\left(\|\phi(\cdot,t,\sigma)\|_X+\|\nabla \phi(\cdot,t)\|_X\right)<\epsilon,
\end{equation}
and
\begin{equation}\label{eq:slowdecay}
\inf_{t\in(-\infty,0)}e^{-\beta t}\|\phi(\cdot,t,\sigma)\|_{L^4(R^3)}>0.
\end{equation}
\end{theorem}

\smallskip
\noindent
{\bf Proof.} Since $\lambda_1(\sigma)$ is differentiable and $\frac{d}{d\sigma}|_{\sigma=\sigma_0}{\rm Re}\,\lambda_1(\sigma)>0 $, we see if we choose $\sigma_1>\sigma_0$ sufficiently close to $\sigma_0$,
$\beta:={\rm Re}\, \lambda_1(\sigma)\in (0,\frac{1}{32})$. By the spectral assumption ({\bf A}), results in Section 2 and spectral mapping theorem (see Section 3.7 page 277 of \cite{Nagel}) we have
\begin{equation}
\sigma(e^{\mathcal{L}_{\sigma_0}})\subseteq B_{e^{-\delta}}\cup \{e^{\lambda_1(\sigma_0)}, e^{\overline{\lambda_1(\sigma_0)}}\}.
\end{equation}
Thus by the continuity of the semigroup $e^{\mathcal{L}_{\sigma}}$ in $\sigma$, we obtain
\begin{equation}
\sigma(e^{\mathcal{L}_{\sigma}})\subseteq B_{e^{-\frac{3}{4}\delta}}\cup \{e^{\lambda_1(\sigma)}, e^{\overline{\lambda_1(\sigma)}}\}
\end{equation}
if $\sigma\in (\sigma_0,\sigma_1)$ and $\sigma_1$ is sufficiently close to $\sigma_0$. Thus
\begin{equation}
\sigma(\mathcal{L}_{\sigma})\subset \{\lambda\in\mathbb{C}:\,{\rm Re}\,\lambda<-\frac{3}{4}\delta\}\cup \{\lambda_1(\sigma),\,\overline{\lambda_1(\sigma)},{\rm~with~Re}\,\lambda_1(\sigma)\in (0,\frac{1}{32})\}.
\end{equation}
Now let us prove the existence of ancient solutions to \ref{eq:ancientequation}. The method is a simple application of the usual techniques of constructing unstable manifold. Our task is easier since we only need to construct one trajectory on the unstable manifold. Let us introduce
\begin{equation}
W:=\{\phi(\cdot,t)\in X, ~t\in(-\infty,0)| \, \esssup_{t\in(-\infty,0)}e^{-\beta t}\left(\|\phi(\cdot,t)\|_X+\|\nabla\phi(\cdot,t)\|_X\right)<\infty\},
\end{equation}
equipped with the natural norm
\begin{equation}
\|\phi\|_{W}:=\esssup_{t\in(-\infty,0)}e^{-\beta t}\left(\|\phi(\cdot,t)\|_X+\|\nabla\phi(\cdot,t)\|_X\right).
\end{equation}
Let us decompose $X=X_u\bigoplus X_s$, the unstable and stable subspaces. Here $X_u$ is two dimensional. Denote $A_u=\mathcal{L}_{\sigma}|_{X_u}$, $A_s=\mathcal{L}_{\sigma}|_{X_s}$; $P_u$ and $P_s$ denotes the projections to the stable and unstable subspaces respectively.
We have the following estimates
\begin{eqnarray}
&&c_1e^{\beta t}\|f\|_X\leq\|e^{A_ut}f\|_X\leq c_2e^{\beta t}\|f\|_X {\rm ~for~}f\in X_u,\, t\in(-\infty,0) {\rm ~and~some~}0<c_1<c_2,\label{eq:unstablepart}\\
&&\|e^{-A_s t}f\|_X\leq C(\sigma)e^{\frac{1}{2}\delta t}\|f\|_X {\rm ~for~}f\in X_s,\, t\in (-\infty,0)\label{eq:stablepart}.
\end{eqnarray}
To find solution to equation (\ref{eq:ancientequation}), for $\phi_{u0}\in X_u$ small, we write the following integral formulation for $t\in(-\infty,0)$
\begin{eqnarray}
&&\phi_u(\cdot,t)=e^{A_u t}\phi_{u0}-\int_0^te^{A_u(t-\tau)}P_u\mathbb{P}(\phi\cdot\nabla\phi)(\cdot,\tau)d\tau,\\
&&\phi_s(\cdot,t)=-\int_{-\infty}^te^{A_s(t-\tau)}P_s\mathbb{P}(\phi\cdot\nabla\phi)(\cdot,\tau)d\tau,{\rm ~~for~}t\in(-\infty,0).
\end{eqnarray}
For $\epsilon>0$ sufficiently small, we will use a contraction mapping argument in $B_{\epsilon}(0)\subset W$ to construct a solution. Firstly, let us assume $\phi_{u0}$ is sufficiently small (note that it's not important to specify the norm, since $X_u$ is two dimensional) such that
\begin{equation}
\|e^{A_u t}\phi_{u0}\|_W<\frac{\epsilon}{2}.
\end{equation}
For $\phi,\,\psi\in W$ let
\begin{eqnarray}
&&T(\phi,\psi)_u(\cdot,t):=\int_0^te^{A_u(t-\tau)}P_u\mathbb{P}(\phi\cdot\nabla\psi)(\cdot,\tau)d\tau,\\
&&T(\phi,\psi)_s(\cdot,t):=\int_{-\infty}^te^{A_s(t-\tau)}P_s\mathbb{P}(\phi\cdot\nabla\psi)(\cdot,\tau)d\tau,{\rm ~~for~}t\in(-\infty,0).
\end{eqnarray}
By the regularity and growth bound of the semigroup, we obtain
\begin{eqnarray*}
&&\|T(\phi,\psi)_u(\cdot,t)\|_X+\|\nabla T(\phi,\psi)_u(\cdot,t)\|_X\\
&&\quad\quad\leq C\int_t^0e^{-\beta(\tau-t)}\max\{1,\frac{1}{\sqrt{\tau-t}}\}e^{2\beta\tau}d\tau \|\phi\|_W\|\psi\|_W\\
&&\quad\quad\quad \leq C(\beta)e^{\beta t}\|\phi\|_W\|\psi\|_W.\\
&&\|T(\phi,\psi)_s(\cdot,t)\|_X+\|\nabla T(\phi,\psi)_s(\cdot,t)\|_X\\
&&\quad\quad\leq C\int_{-\infty}^te^{-\frac{1}{2}\delta (t-\tau)}\max\{1,\sqrt{t-\tau}\}e^{2\beta\tau}d\tau\|\phi\|_W\|\psi\|_W\\
&&\quad\quad\quad\leq C(\delta,\beta)e^{2\beta t} \|\phi\|_W\|\psi\|_W.
\end{eqnarray*}
Thus we conclude
\begin{equation}
\|T(\phi,\psi)\|_W\leq C(\delta,\beta)\|\phi\|_W\|\psi\|_W.
\end{equation}
Now it's routine to check if we choose $\epsilon$ sufficiently small, we can use contraction mapping argument in $B_{\epsilon}\subset W$ to obtain a solution $\phi(\cdot,t,\sigma)\in B_{\epsilon}\subset W$. We only need to verify equation (\ref{eq:slowdecay}). We can make $\|\phi_{u0}\|_X=\gamma \epsilon$ for some small positive $\gamma$. Note
\begin{eqnarray*}
&&\|\phi_u(\cdot,t)-e^{A_u t}\phi_{u0}\|_X+\|\nabla(\phi_u(\cdot,t)-e^{A_u t}\phi_{u0})\|_X\leq C(\beta)e^{\beta t} \epsilon^2,\\
&&\|\phi_s(\cdot,t)\|_X+\|\nabla\phi_s(\cdot,t)\|_X\leq C(\beta,\delta)e^{\beta t}\epsilon^2, {\rm ~for~}t\in(-\infty,0).
\end{eqnarray*}
Thus
\begin{eqnarray*}
&&\|\phi(\cdot,t)\|_{L^4(R^3)} \ge \|\phi_u(\cdot,t)\|_{L^4(R^3)}-\|\phi_s(\cdot,t)\|_{L^4(R^3)}\\
&&\quad\quad\ge\|e^{A_u t}\phi_{u0}\|_{L^4(R^3)}-C(\beta,\delta)e^{\beta t}\epsilon^2\ge c(\alpha,\delta,\gamma)\epsilon e^{\beta t}, {\rm ~for~}t\in(-\infty,0),
\end{eqnarray*}
if we choose $\epsilon$ sufficiently small. Thus the theorem is proved.\\

Now let us consider the situation that the spectral assumption ({\bf B}) holds. In this case, we have the following theorem.
\begin{theorem}\label{th:bifurcationB}
Suppose the spectral assumption ({\bf B}) holds. Let $\mathcal{L}_{\sigma_0}$, $U_{\sigma_0}$ and $U$ be defined as above. Denote $\mathsf{v}\in\mathcal{D}$ as a unit eigenfunction for $\mathcal{L}_{\sigma_0}$ corresponding to eigenvalue $0$. In addition, assume the following generic non-degeneracy conditions
\begin{eqnarray}
&&\mathbb{P}(U_{\sigma_0}\cdot\nabla U+U\cdot\nabla U_{\sigma_0})\notin {\rm Range}\,(\mathcal{L}_{\sigma_0}),\label{eq:nondegeneracy1}\\
&&\mathbb{P}(\mathsf{v}\cdot\nabla \mathsf{v})\notin {\rm Range}\,(\mathcal{L}_{\sigma_0})\label{eq:nondegeneracy2}.
\end{eqnarray}
Then there exists a solution curve $(U_{\sigma(s)},\sigma(s))$ to equations (\ref{eq:forwardselfsimilarequation1},\ref{eq:forwardselfsimilarequation2}), for $s\in(-\epsilon,\epsilon)$ with some $\epsilon>0$. Moreover, $\sigma(0)=\sigma_0, ~\frac{d}{ds}\sigma(0)=0$ and $\sigma(s)$ is twice continuously differentiable with $\frac{d^2}{ds^2}\sigma(0)<0$. Thus for $\sigma<\sigma_0$ sufficiently close to $\sigma_0$ we have two forward self similar solutions to Navier Stokes equation with initial data $\sigma u_0$.
\end{theorem}

\noindent
{\bf Remark.} Since ${\rm Range}\,(\mathcal{L}_{\sigma_0})$ is closed with co-dimension $1$, the assumptions in the theorem are generic.

\smallskip
\noindent
{\bf Proof.} Again to prove the theorem, we seek solution $U_{\sigma}$ in the form
\begin{equation}
U_{\sigma}=U_{\sigma_0}+(\sigma-\sigma_0)U+\phi, {\rm ~with~}\phi\in X.
\end{equation}
Then $\phi$ satisfies $F(\phi,\sigma)=0$, where
\begin{eqnarray*}
&F(\phi,\sigma):=\Delta \phi+\frac{x}{2}\cdot\nabla \phi+\frac{\phi}{2}-U_{\sigma_0}\cdot\nabla\phi-\phi\cdot\nabla U_{\sigma_0}-(\sigma-\sigma_0)U\cdot\nabla \phi-(\sigma-\sigma_0)\phi\cdot\nabla U\\
&\quad\quad -(\sigma-\sigma_0)U\cdot\nabla U_{\sigma_0}-(\sigma-\sigma_0)U_{\sigma_0}\cdot\nabla U-(\sigma-\sigma_0)^2U\cdot\nabla U-\phi\cdot\nabla\phi+\nabla P.
\end{eqnarray*}
It is easy to verify that $F:\,\mathcal{D}\times (\sigma_0-\epsilon,\sigma_0+\epsilon)\rightarrow X$ is smooth, where $\epsilon>0$ is a small number. In addition we can verify
\begin{eqnarray}
&&\partial_{\sigma}F(0,\sigma_0)\notin {\rm Range}\,(\partial_{\phi}F(0,\sigma_0)),\\
&&\partial^2_{\phi\phi}F(0,\sigma_0)[\mathsf{v},\mathsf{v}]\notin {\rm Range}\,(\partial_{\phi}F(0,\sigma_0)).
\end{eqnarray}
Thus we can apply the saddle-node bifurcation theorem (see Theorem I.4.1 page 12 from \cite{Hans}) and finish the proof. \\

The above two cases with spectral assumptions ({\bf A}) or ({\bf B}) are generic if we indeed have eigenvalues of $\mathcal{L}_{\sigma}$ crossing the imaginary lines when we increase $\sigma$. There is an important non-generic case when the initial data $u_0$ is axi-symmetric in addition to being $-1$ homogeneous. We will not give the full details, but just sketch some of the ideas. Suppose when we increase $\sigma$ to $\sigma_0$ the spectrum of $\mathcal{L}_{\sigma}$ splits into two parts, one strictly to the left of imaginary axis, the other consists of an eigenvalue crossing $0$. Now suppose for $\sigma=\sigma_0$ the eigenfunctions corresponding to $0$ eigenvalue of $\mathcal{L}_{\sigma_0}$ are not symmetric (note the solution $U_{\sigma}$ for $\sigma<\sigma_0$ is axi-symmetric). Then since the restriction of $\mathcal{L}_{\sigma_0}$ to axi-symmetric vector fields is invertible, we can continue the axi-symmetric solution curve $U_{\sigma}$ in a regular way. In addition, we can expect a ``pitchfork" bifurcation, which gives additional solutions. Thus we have non-uniqueness of forward self similar solutions with initial data $\sigma>\sigma_0$ sufficiently close to $\sigma_0$ in this case as well.

\end{section}

\begin{section}{Localize solutions with $-1$ homogeneous initial data}
In this section we will localize the different solutions to Navier Stokes equation with $-1$ homogeneous initial data $\sigma u_0$, which appear through bifurcations under spectral assumption ({\bf A}), or ({\bf B}) together with the non-degeneracy conditions in Theorem \ref{th:bifurcationB}. We will only give proof in the case of spectral assumption ({\bf A}) below, as the other case is identical.
Under the spectral condition ({\bf A}), we see for $\sigma_0<\sigma<\sigma_1$ there are at least two solutions to NSE with scale-invariant initial data $\sigma u_0\in C^{\infty}(R^3\backslash\{0\})$, one scale invariant
and is given by
\begin{equation}
u_1(x,t)=\frac{1}{\sqrt{t}}U_{\sigma}(\frac{x}{\sqrt{t}}).
\end{equation}
The other solution is not self similar and is given by
\begin{equation}
u_2(x,t)=\frac{1}{\sqrt{t}}U_{\sigma}(\frac{x}{\sqrt{t}})+\frac{1}{\sqrt{t}}\phi(\frac{x}{\sqrt{t}},\log{t}),
\end{equation}
where $\phi(\cdot,s)$ is an ancient solution to equation (\ref{eq:ancientequation}) and satisfies equations (\ref{eq:smallancientsolution},\ref{eq:slowdecay}). If we take $\sigma$ sufficiently close to $\sigma_0$, we can assume $s(\mathcal{L}-\mathcal{K}(U_{\sigma}))\leq \frac{1}{32}$. Moreover, we can assume $\sup_{s\in (-\infty,0)}(\|\phi(\cdot,s)\|_X+\|\nabla\phi(\cdot,s)\|_X)$ is sufficiently small so that we can apply Theorem \ref{th:perturbingsingularsolution} with $a=U_{\sigma}$ and $b=\phi$ again by making $\sigma$ sufficiently close to $\sigma_0$. Now fix such a $\sigma$. Let us decompose
\begin{equation}
u=v_0+w_0,
\end{equation}
where $v_0$ is compactly supported, divergence free and $v_0|_{B_R(0)}\equiv u_0|_{B_R(0)}$. We can make $\|w_0\|_{L^4(R^3)}\leq C(\sigma_0,u_0)R^{-\frac{1}{4}}$.
Now apply Theorem \ref{th:perturbingsingularsolution} first with initial data $-w_0$, $a=U_{\sigma}$ and $b=0$. By taking $R$ sufficiently large, we obtain a solution $\tilde{u}_1(x,t)$ to equation (\ref{eq:singularstokes}) in $R^3\times(0,1)$. One can easily check $u_1+\tilde{u}_1$ is a solution to Navier Stokes equation with initial data $v_0$. Similarly, apply Theorem \ref{th:perturbingsingularsolution} with initial data $-w_0$, $a=U_{\sigma}$ and $b=\phi$, we obtain a solution $\tilde{u}_2$ to equation (\ref{eq:singularstokes}) in $R^3\times (0,1)$ (taking $R$ appropriately in the first step). Thus $u_2+\tilde{u}_2$ is a solution to Navier Stokes equation with initial data $v_0$. In fact by the estimates on $u_1,~\tilde{u}_1,~u_2,~\tilde{u}_2$ and the regularity theory of NSE, we can easily deduce that $u_1+\tilde{u}_1,\,u_2+\tilde{u}_2\in C^{\infty}(R^3\times (0,1))$. Using $v_0\in L^2(R^3)$, it's not hard to obtain that the two solutions are uniformly bounded in $L^2(R^3)$ in time. Lastly we shall show $u_1+\tilde{u}_1$ is not identically equal to $u_2+\tilde{u}_2$. Since $\tilde{u}_1(\cdot,t)$ and $\tilde{u}_2(\cdot,t)$ are bounded in $L^4(R^3)$. We see
\begin{equation}
\|(u_1+\tilde{u}_1-u_2-\tilde{u}_2)(\cdot,t)\|_{L^4(R^3)}\ge \|(u_1-u_2)(\cdot,t)\|_{L^4(R^3)}-C\ge t^{-\frac{1}{8}}\|\phi(\cdot,\log{t})\|_{L^4(R^3)}-C
\end{equation}
is unbounded as $t\to 0+$ by (\ref{eq:slowdecay}). Thus we have proved the following theorem.
\begin{theorem}
Assume the spectral condition ({\bf A}) holds. Then there exist two different Leray-Hopf weak solutions which are smooth in $R^3\times(0,1)$ with the same initial data $v_0\in C^{\infty}(R^3\backslash\{0\})$. $v_0(x)=O(\frac{1}{|x|})$ near origin.
\end{theorem}

Similarly, we can obtain the following theorem.
\begin{theorem}
Assume the spectral condition ({\bf B}) and the non-degeneracy conditions (\ref{eq:nondegeneracy1},\ref{eq:nondegeneracy2}) hold. Then there exists compactly supported divergence free
$v_0\in C^{\infty}(R^3\backslash\{0\})$ with $v_0(x)=O(\frac{1}{|x|})$ near origin, such that there are two Leray-Hopf weak solutions with initial data $v_0$. Moreover the two Leray-Hopf solutions are smooth in $R^3\times (0,\infty)$.
\end{theorem}

\end{section}

\end{document}